\documentclass[12pt]{article}
\def\hind{\hangindent=2pc\hangafter=1}
\newfont{\smcaps}{cmcsc10 scaled\magstep1}

\newcommand{\MA}{{\rm ~MA\,}}
\newcommand{\ARMA}{{\rm ~ARMA\,}}
\newcommand{\AR}{{\rm ~AR\,}}

\font\tty=cmtt10 at 11truept
\usepackage{epsfig}
\raggedright
\begin{document}
\baselineskip=22pt

\title{Improved Pe\v{n}a-Rodriguez Portmanteau Test}
\date{June 19, 2006}
\author{{Jen-Wen Lin$^{\rm a\/}$ and A. Ian McLeod$^{\rm a\/}$}\\\\
$^{\rm a\/}$\it Department of Statistical and Actuarial Sciences,\\
\it The University of Western Ontario,\\
\it London, Ontario Canada N6A 5B7\\
}
\maketitle
\hrule
\newpage
{\bf Abstract\/}

\bigskip
\noindent
Several problems with the diagnostic check suggested by Pe\v{n}a and Rodriguez (2002)
are noted and an improved Monte-Carlo version of this test is suggested.
It is shown that quite often
the test statistic recommended by Pe\v{n}a and Rodriguez (2002) may not exist
and their asymptotic distribution of the test does not agree with the suggested
gamma approximation very well
if the number of lags used by the test is small.
It is shown that the convergence of this test statistic to its asymptotic distribution
may be quite slow when the series length is less than 1000
and so a Monte-Carlo test is recommended.
Simulation experiments suggest the Monte-Carlo test is usually more powerful than
the test given by Pe\v{n}a and Rodriguez (2002) and often much more powerful
that the Box-Ljung portmanteau test.
Two illustrative examples of enhanced diagnostic checking with the Monte-Carlo test
are given.

\bigskip
\hrule

\newpage
{\noindent \bf 1. Introduction}
\medskip

Let $X_t, t=1,2,\ldots$ be a stationary and invertible {\ARMA($p,q$)} model (Box, Jenkins, and Reinsel, 1994),
\begin{equation}
(1-\phi_{1}B-\cdots -\phi_{p}B^{p}) X_{t}=( 1-\theta_{1}B-\cdots -\theta_{q}B^{q}) a_{t},
\newcounter{ARMA}
\setcounter{ARMA}{\value{equation}}
\end{equation}
where $B$ is the backshift operator on $t$, $a_t$ is a sequence of independent and identical
normal random variables with mean zero and variance $\sigma_a^2$.
After fitting this model to a series of length $n$, the residual autocorrelations,
\begin{equation}
\hat{r}(k)=\sum_{t=k+1}^{n} \hat{a}_{t}\hat{a}_{t-k}/\sum_{t=1}^{n}\hat{a}_{t}^{2}\ \ \ \ k=1,2,\cdots,
\newcounter{racf}
\setcounter{racf}{\value{equation}}
\end{equation}
where $\hat{a}_t$ denotes the fitted residuals, may be used
for checking model adequacy.
One of the most widely used model diagnostic checks (Li, 2004) is the portmanteau test of
Ljung and Box (1978),
\begin{equation}
Q_m = n (n+2) \sum_{k=1}^{m} (n-k)^{-1} \hat r(k)^{2},
\newcounter{LjungBox}
\setcounter{LjungBox}{\value{equation}}
\end{equation}
where, under the assumption of model adequacy, $Q_m$ is approximately
$\chi^2_{m-(p+q)}$ distributed.

A new portmanteau diagnostic test statistic,
$\hat D_m = n(1- |\hat{R}_{m}|^{1 / m})$,
based on the determinant of the residual autocorrelation matrix,
\begin{equation}
\hat{R}_{m}=\left( \begin{array}{cccc}
 1 & \hat{r}(1) & \cdots  & \hat{r}(m) \\
 \hat{r}(1) & 1 & \cdots  & \hat{r}(m-1) \\
 \vdots  & \cdots  & \ddots  & \vdots  \\
 \hat{r}(m) & \cdots  & \hat{r}(1) & 1
\end{array}\right)  ,
\end{equation}
was suggested by Pe\v{n}a and Rodriguez (2002).
As noted by Pe\v{n}a and Rodriguez (2002),
$|\hat{R}_{m}|$ is the estimated generalized
variance of the residuals standardized by dividing by their standard deviation.

\medskip
{\noindent \bf 2. THE $D_m$ TEST AND ITS LIMITATIONS}
\medskip

Pe\v{n}a and Rodriguez (2002, Theorem 1) showed that
if the model is correctly identified,
$\hat{D}_m$ is asymptotically distributed as
$\sum \limits_{i=1}^{m}\lambda_{i}\chi _{1,i}^{2}$,
where $\chi _{1,i}^{2}$ are independent
Chi-squared random variables with one degree of freedom,
and $\lambda _{i}$ are the eigenvalues of ${\cal Q}_m W_{m}$, where
$W_{m}$ is a diagonal matrix with the $i$-th diagonal elements, $w_{i}=(m-i+1)/m$,
$i=1,2,\ldots, m$ and ${\cal Q}_m$ is the asymptotic covariance matrix of the normalized
residual autocorrelations $\surd{n} (\hat{r}(1), \ldots, \hat{r}(m))$ given
in McLeod (1978, eqn. 15).
As pointed out by Pe\v{n}a and Rodriguez (2002), this asymptotic distribution
may be computed using the method of Imhof (1961).
We have implemented the computation of this asymptotic distribution for
the general ARMA$(p,q)$ model in our R package {\tty gvtest}.
The cumulative distribution function for this asymptotic distribution
may be denoted by $F(x;\lambda_1,\ldots,\lambda_m)$.

Pe\v{n}a and Rodriguez (2002) suggested evaluating $\hat{D}_m$
by a gamma approximation distribution to the asymptotic distribution.
The approximation distribution is derived
by equating the first two moments of a gamma distribution with those
of the corresponding asymptotic distribution.
The density function for this gamma approximation may be written
$f_\gamma(x;\alpha,\beta)=
e^{-\frac{x}{\beta }} x^{\alpha -1} \beta ^{-\alpha } / \Gamma (\alpha)$
where,
\begin{equation}
\alpha =3\left. m \{\left( m+1\right) -2\left( p+q\right) \right\}
^{2} /\left[ 2\{ 2\left( m+1\right) \left( 2m+1\right) -12m( p+q)
\} \right],
\end{equation}
and
\begin{equation}
\beta =3\left. m \{\left( m+1\right) -2\left( p+q\right) \right\}
/\left\{ 2\left( m+1\right) \left( 2m+1\right) -12m( p+q) \right\}.
\end{equation}
Pe\v{n}a and Rodriguez (2002) indicated that the approximation
improves as $m$ increases.
In addition, it may be noted that if $m$ is too small then it can happen
that $\alpha \le 0$ or $\beta \le 0$ which is numerically infeasible.
For example if $p+q=3$ then we must have $m \ge 8$ to make $\alpha>0$
and $\beta>0$.

Pe\v{n}a and Rodriguez (2002) found that the empirical distribution
of $\hat D_m$ did not agree very well with the gamma approximation so they
suggested a modified statistic,
\begin{equation}
D_m=n-n |\ddot {R}_{m} | ^{\frac{1}{m}},
\end{equation}
where $\ddot{R}_{m}$ denotes the residual autocorrelation matrix
replacing $\hat{r}(k)^2$ with $\ddot{r}_{k}^{2}$, where
\[
\ddot{r}(k)^{2}=\left( n+2\right) \left( n-k\right) ^{-1} \hat{r}(k)^{2}.
\]
This is similar in spirit to the modification suggested by Ljung and Box (1978)
to the original Box and Pierce portmanteau statistic (1970).
Although, as shown by Pe\v{n}a and Rodriguez (2002),
this approximation works well when $n \le 100$
in first order autoregressive models,
it does not provide a good approximation to the asymptotic distribution
for more complicated models if the number of lags, $m$, is small.
As shown in Figure 1, the gamma approximation can distort the size of a 5\% significance test
relative to the asymptotic distribution for {\AR}$(2)$ models when $m$ is $10$.
Similar distortions are found for {\MA}$(2)$ models
or higher orders {\AR} and {\MA} models.
The distortions for {\ARMA}$(1,1)$ with $m=10$ models were also investigated and
the results were listed in Table 1.
Moreover, this distortion would tend to make the tests
based on the gamma approximation reject more often than they should.
In other words,
the test based on the gamma approximation is not conservative.
So despite the fact that as shown by Pe\v{n}a and Rodriguez (2002, Table 2) the small
sample performance is acceptable in some cases, the more general use
of tests based on the gamma approximation can not be recommended.
The reader may investigate the accuracy of the gamma approximation for a
particular $m$ and ARMA$(p,q)$  using our {\tty gvtest} package.

\bigskip
\begin{center}
[Table 1 about here]
\end{center}
\medskip

Another serious limitation to the use of $D_m$ is that it is frequently undefined.
This happens because $\ddot{R}_{m}$ is not always positive-definite (McLeod and Jim\' enez, 1984).
When we tried to replicate Table 4 in Pe\v{n}a and Rodriguez (2002) we frequently found cases
where $\ddot{R}_{m}$ was not positive-definite.
For example, for model 7 listed in this table this happens about 25\% of the time with $m=10$.
This problem occurred with many other models listed in this table.
Although only very short time series were used in Table 4, this problem also occurs
with longer time series particularly when $m$ is also larger.
For this reason, it is better to concentrate on the original $\hat D_m$ statistic.

\bigskip
\begin{center}
[Figure 1 about here]
\end{center}
\medskip

{\noindent \bf 3. TESTS BASED ON $\hat D_m$}
\medskip

In a modern high level computing environment, such as R or {\it Mathematica\/},
it is not difficult to evaluate
the asymptotic distribution and so a test based directly on this asymptotic distribution might seem
to be preferable.
Unfortunately, the convergence of $\hat D_m$ to its asymptotic distribution is often very slow.
In Table 2, we evaluated the asymptotic distribution corresponding to the upper 5\% point
of the empirical distribution of $\hat D_m$ for the first-order autoregressions
with series lengths $n=100, 200, 500$ and 1000.
Not until $n$ is very large is the asymptotic distribution reliable.
Figure 2 shows a QQ plot for these simulations which shows that the discrepancy between the
empirical and asymptotic distribution increases at the larger quantiles.
The asymptotic distribution tends to understate the actual finite-sample significance
level while the gamma approximation errs in the opposite direction.

For these reasons, a Monte-Carlo test procedure (Gentle, 2002, \S 2.3) is recommended when $n<1000$.
This procedure is quite practical on typical computers now available.
This test is essentially equivalent to a parametric bootstrap test (Davison and Hinkley, 1997, Ch.4).
The steps in the procedure are indicated below:

\begin{enumerate}
\item After fitting model in eqn. (\theARMA) obtain $\hat D_m$.
\item Select the number of Monte-Carlo simulations, $N$. Typically $100 \le N \le 1000$.
\item Simulate the model in eqn. (\theARMA) using the estimated parameters obtained
in Step 1 and obtain $\hat D_m$ after estimating the parameters in the simulated series.
\item Repeat Step 3 $N$ times counting the number of times $k$ that a value
of $\hat D_m$ greater than or equal to that in Step 1 has been obtained.
\item The $P$-value for the test is $(k+1)/(N+1)$.
\item Reject the null hypothesis if the $P$-value is smaller than a predetermined significance level.
\end{enumerate}

It should be noted that nuisance parameters are present in our proposed procedure
and this could cause size distortion in the Monte Carlo test (J\"ockel, 1986).
The empirical size of the Monte-Carlo test for the first order autoregressive model
was investigated by simulation.
The results were summarized in Table 3 and it is seen that
the empirical sizes are very close to their nominal level.
In general, it has been shown
(Dufour, 2006; Dufour and Khalaf, 2001) that if consistent estimators
are used then the Monte-Carlo test produces an asymptotically correct size as $n \rightarrow \infty$.
The asymptotic validity of the Monte-Carlo test follows immediately.
Alternatively since the gamma approximation is asymptotically correct as both
$n$ and $m$ get large, it follows that $\hat D_m$ is asymptotically pivotal and
hence does not depend on nuisance parameters.

\bigskip
\begin{center}
[Figure 2 and Table 2,3 about here]
\end{center}

\medskip
{\noindent \bf 4. EMPIRICAL POWER COMPARISONS}
\medskip

In many cases the Monte-Carlo test outperforms the the $D_m$ test
based on the gamma approximation.
As an illustration, in Table 4 the simulation results for the Monte-Carlo test
are compared with $D_m$ test for the GARCH models Pe\v{n}a and Rodriguez (2002, Table 12).
We see that the Monte-Carlo test has always has higher power than the gamma approximation
in Table 4.

Pe\v{n}a and Rodriguez (2002) indicated that
the advantage of their test over the Ljung-Box test, denoted as $Q_m$,
may disappear in heteroscedastic data with long persistence.
The results in Table 5 confirm this fact.
For example,
$\hat D_m$ has power $0.888$ for $m=5$ and $d=0.3$ with respect to a
power of $0.880$ for $Q_m$.
However, it is interesting to note that, as can be seen in Table 5,
the difference in power increases as $n$ and $m$ increase, so that,
for example, when $n=512, m=40$, $\hat D_m$ is about $100\times (0.446-0.377)/0.377 \doteq 18\%$ more powerful
than $Q_m$.

In another simulation experiment, shown in Table 6,
we compared the power of Monte-Carlo tests using both
$\hat D_m$ and $Q_m$ for twelve models examined by
Pe\v{n}a and Rodriguez (2002, Table 3).
Overall our results are similar to those reported by Pe\v{n}a and Rodriguez (2002, Table 3).
There are some differences though and this is due the limitations discussed in \S 2.
Specifically, incorrect size when the gamma approximation was used or bias in the
simulations caused by the fact that $D_m$ frequently does not exist.
We also investigated Monte-Carlo tests based on the Ljung-Box test, $Q_m$.
As shown in Table 6, the Monte-Carlo test using $\hat D_m$ outperforms this test.
The empirical power for the Monte-Carlo test for the Box-Pierce test (Box and Pierce, 1970)
was also computed but it was not significantly different from the results for the
$Q_m$ test.

\bigskip
\begin{center}
[Table 4, 5 and 6 about here]
\end{center}

\medskip
{\noindent \bf 5. ILLUSTRATIVE EXAMPLES}
\medskip

[Note: {\tty GVTest} is currently available on our website:
http://www.stats.uwo.ca/faculty/aim/2005/GVTest/
and it will be put on CRAN when our paper is published.]

Our Monte-Carlo test is implemented in an R package
available on CRAN, {\tty GVTest}.
Hipel and McLeod (1978) fit an ARMA$(2,1)$ model a tree-ring time series denoted by
{\tty Ninemile}.  There were $n=771$ annual values.
As can be seen from Table 7,
if one uses the Ljung-Box portmanteau test with $m=40$, the model appears
adequate although for $m=20$ the test does strongly suggest model inadequacy.
In this case the $\hat D_m$ Monte-Carlo test provides a clearer indication since
it indicates to reject for $m=20, 30, 40, 50$.

Fitting an AR(2) model to the sunspot series in Box, Jenkins and Reinsel (1994, Series E) we
again found that the Ljung-Box test suggests the model is adequate but the Monte-Carlo $\hat D_m$
test indicates model inadequacy.
Note that with $m=5$, the Ljung-Box test does have a $p$-value of about 5\% but usually $m$ is taken
to be larger since only for large enough $m$ is the covariance matrix idempotent.
For this reason, the result for $m=5$ might be discounted.
When a Monte-Carlo test is used, no such difficulties arise.
The results are summarized in Table 7.

\bigskip
\begin{center}
[Table 7 about here]
\end{center}

\medskip
{\noindent \bf 8. CONCLUDING REMARKS}
\medskip

The implementation of the generalized variance
portmanteau test statistic suggested by Pe\v{n}a and Rodriguez (2002) is unsatisfactory because
it frequently does not exist and there are important limitations to the gamma approximation.
These difficulties are rectified by using a Monte-Carlo test.
Further simulation experiments have indicated that the Monte-Carlo portmanteau test using
simulated Gaussian innovations often works well even when the true is non-normal.
This is the case for thicker tail distributions such as double exponential and
the $t$ distribution on 5 df.

\medskip
{\noindent ACKNOWLEDGEMENT}
\medskip

A.I. McLeod acknowledges with thanks a Discovery Grant Award from NSERC.

\newpage
\centerline{REFERENCES}
\parindent 0pt
\medskip

\hind
{Box, G.E.P., Jenkins, G.M. and Reinsel, G.C.\/}  (1994),
{\it Time Series Analysis: Forecasting and Control\/}, 3rd Ed.,
San Francisco: Holden-Day.

\hind
{Box, G.E.P. and Pierce, D.A.} (1970),
``Distribution of the residual autocorrelation in autoregressive integrated moving average time series models,''
{\it Journal of American Statistical Association\/}, {65}, 1509--1526.

\hind
Davison, A.C. and Hinkley, D.V. (1997),
{\it Bootstrap Methods and their Application},
Cambridge: Cambridge University Press.

\hind
Dufour, J.-M. (2006, to appear).
``Monte Carlo Tests with Nuisance Parameters :
A General Approach to Finite-Sample Inference and Nonstandard Asymptotics in Econometrics'',
{\it Journal of Econometrics}.

\hind
{Dufour, J.-M. and Khalaf, L.} (2001),
``Monte-Carlo Test Methods in Econometrics,''
In {\it Companion to Theoretical Econometrics\/},
Ch. 23, 494--519.
Oxford: Blackwell.

\hind
Gentle, J.E. (2002),
{\it Elements of Computational Statistics,}
New York: Springer.

\hind
Hipel, K.W. and McLeod, A.I. (1978),
``Preservation of the rescaled adjusted range, Part 2, Simulation studies using Box-Jenkins models,''
{\it Water Resources Research}, {14}, 509--516.

\hind
{Imhof, J.P.} (1961),
``Computing the distribution of quadratic forms in normal variables,''
{\it Biometrika }, {48}, 419--426.

\hind
{J\"ockel, K.F.} (1986),
``Finite sample properties and asymptotic efficiency of Monte Carlo tests,''
{\it Annals of Statistics\/}, 14, 336--347.

\hind
{Li, W.K.} (2004),
{\it Diagnostic Checks in Time Series\/},
New York: Chapman and Hall/CRC.

\hind
{Ljung, G.M.}  (1986),
``Diagnostic testing of univariate time series models,''
{\it Biometrika\/}, {73}, 725--730.

\hind
{Ljung, G.M. and Box, G.E.P.} (1978).
``On a Measure of Lack of Fit in Time Series Models,''
{\it Biometrika\/}, {65}, 297--303.

\hind
{McLeod, A.I.} (1978),
``On the distribution of residual autocorrelations in Box-Jenkins models,''
{\it Journal of the Royal Statistical Society B\/}, {40}, 396--402.

\hind
{McLeod, A.I. and Jim\' enez, C.} (1984),
``Nonnegative definiteness of the sample autocorrelation function,''
{\it The American Statistician}, {38}, 297--298.

\hind
{Pe\v{n}a, D. and Rodriguez, J.} (2002),
``A powerful portmanteau test of lack of fit for time series,''
{\it Journal of American Statistical Association\/}, {97}, 601--610.

\newpage

Table 1.
P-values of the upper 5 $\%$ quantiles of the gamma approximation
for {\ARMA}$(1,1)$ models with
$\phi_1$ and $\theta_1=\pm 0.9,\pm 0.6, \pm 0.3$ evaluated by
their asymptotic distributions.

\begin{center}
\begin{tabular}{ccccccc}
\noalign{\smallskip}
\noalign{\hrule}
\noalign{\smallskip}
\noalign{\hrule}
\noalign{\smallskip}
$\theta_1|\ \phi_1$&$-0.9$&$-0.6$&$-0.3$&$0.3$&$0.6$&$0.9$ \\
\noalign{\smallskip}
\noalign{\hrule}
\noalign{\smallskip}
$-0.9$&     &$0.105$&$0.091$&$0.083$&$0.085$ &$0.109$ \\
$-0.6$&$0.105$&     &$0.069$&$0.063$&$0.065$ &$0.085$ \\
$-0.3$&$0.091$&$0.692$&   &$0.060$&$0.063$ &$0.083$ \\
$0.3$ &$0.083$&$0.063$&$0.060$&   &$0.069$ & $0.091$ \\
$0.6$ &$0.085$&$0.065$&$0.063$&$0.069$&    & $0.105$ \\
$0.9$ &$0.108$&$0.085$&$0.083$&$0.091$&$0.105$ &  \\
\noalign{\smallskip}
\noalign{\hrule}
\end{tabular}
\end{center}

\newpage

Table 2:
Asymptotic probability corresponding to the upper 5\% quantile of the
empirical distribution of $\hat D_m$ with $m=50$ for the first-order
autoregressive model with parameter $\phi$ and series length $n$.
Each entry in the table is based on $10^4$ simulations.

\begin{center}
\begin{tabular}{ccccc}
\noalign{\hrule}
\noalign{\smallskip}
\noalign{\hrule}
\noalign{\smallskip}
&\multicolumn{4}{c}{$\phi$} \\
\noalign{\smallskip}
\cline{2-5}
\noalign{\smallskip}
$n$&$-0.8$&$-0.4$&$0.4$&$0.8$ \\
\noalign{\smallskip}
\noalign{\hrule}
\noalign{\smallskip}
$100$&$0.423$&$0.433$&$0.439$&$0.445$ \\
$200$&$0.179$&$0.161$&$0.170$&$0.181$ \\
$500$&$0.093$&$0.092$&$0.081$&$0.079$ \\
$1000$&$0.060$&$0.059$&$0.051$&$0.062$ \\
\noalign{\smallskip}
\noalign{\hrule}
\end{tabular}
\end{center}

\newpage
Table 3:
Empirical significance levels of $\hat D_m$ under first order autoregressive models.
The empirical power for the MC test is based on $1,000$ simulations.
Each MC test also used $250$ simulations.

\begin{center}
\begin{tabular}{cccccc}
\noalign{\smallskip}
\noalign{\hrule}
\noalign{\smallskip}
\noalign{\hrule}
\noalign{\smallskip}
&\multicolumn{4}{c}{$\alpha=0.05$} \\
\noalign{\smallskip}
%\cline{2-6}
%\noalign{\smallskip}
$             $ & $\phi=0.1$ &$\phi=0.3$  &$\phi=0.5$     &$\phi=0.7$   &  $\phi=0.9$ \\
$m=10$        & $0.050$    &$0.049$     &$0.059$        &$0.055$      &  $0.049$ \\
$m=20$        & $0.049$    &$0.047$     &$0.046$        &$0.050$      &  $0.058$  \\
\noalign{\smallskip}
\noalign{\hrule}
\noalign{\smallskip}
\noalign{\smallskip}
&\multicolumn{4}{c}{$\alpha=0.01$} \\
\noalign{\smallskip}
%\cline{2-6}
%\noalign{\smallskip}
$           $ & $\phi=0.1$ &$\phi=0.3$  &$\phi=0.5$     &$\phi=0.7$   &  $\phi=0.9$ \\
$m=10$        & $0.016$    &$0.014$     &$0.011$        &$0.012$      &  $0.008$ \\
$m=20$        & $0.009$    &$0.012$     &$0.014$        &$0.014$      &  $0.012$  \\
\noalign{\smallskip}
\noalign{\hrule}
\end{tabular}
\end{center}

\newpage
Table 4:
Power Comparison of MC Test and $D_m$ Test for GARCH Time Series.  The empirical power for the MC test is based on $1,000$ simulations.
Each MC test also used $1,000$ simulations.  The empirical power reported by Pe\v{n}a and Rodriguez (2002, Table 9) is
shown in the column $D_m$.  Models $A$ and $B$ refer to the two GARCH models used by Pe\v{n}a and Rodriguez (2002).

\begin{center}
\begin{tabular}{ccccccccccc}
\hline\hline
&&&\multicolumn{2}{c}{\rule[-3mm]{0mm}{8mm}$m=12$}&&\multicolumn{2}{c}{$m=24$}&&\multicolumn{2}{c}{$m=32$}\\
\noalign{\smallskip}
\cline{4-5}\cline{7-8}\cline{10-11}
\noalign{\smallskip}
Model&$n$&\ &$\hat D_m$&$D_m$&\ &$\hat D_m$&$D_m$&\ &$\hat D_m$&$D_m$\\
\noalign{\smallskip}
%\noalign{\hrule}
%\noalign{\smallskip}
A&$250$&\ &$0.387$&$0.268$&\ &$0.397$&$0.244$&\ &$0.381$&$0.213$\\
B&$250$&\ &$0.878$&$0.821$&\ &$0.859$&$0.782$&\ &$0.841$&$0.731$\\
A&$500$&\ &$0.574$&$0.522$&\ &$0.566$&$0.501$&\ &$0.566$&$0.479$\\
B&$500$&\ &$0.992$&$0.986$&\ &$0.990$&$0.981$&\ &$0.987$&$0.973$\\
A&$1000$&\ &$0.843$&$0.807$&\ &$0.843$&$0.802$&\ &$0.829$&$0.776$\\
B&$1000$&\ &$1.000$&$1.000$&\ &$1.000$&$1.000$&\ &$1.000$&$1.000$\\
\noalign{\smallskip}
\noalign{\hrule}
\end{tabular}
\end{center}

\newpage

Table 5.
Empirical power comparison of Monte-Carlo test using $\hat D_m$ and Ljung-Box test, $Q_m$,
for fractional noise time series
with $n=256, 512$ and $d= 0.2, 0.3$.  The empirical power is based on $1,000$ simulations and each
Monte-Carlo test also uses $1,000$ replications.
The first entry corresponds to $\hat D_m$ and the second to $Q_m$.

\begin{center}
\begin{tabular}{ccccccc}
\noalign{\smallskip}
\noalign{\hrule}
\noalign{\smallskip}
\noalign{\hrule}
\noalign{\smallskip}
$d$&$n$&$m=5$&$m=10$&$m=20$&$m=30$&$m=40$ \\
\noalign{\smallskip}
\noalign{\hrule}
\noalign{\smallskip}
$0.2$&$256$&$0.283/0.274$&$0.274/0.247$&$0.227/0.202$&$0.192/0.176$&$0.171/0.153$ \\
$0.3$&$256$&$0.539/0.517$&$0.540/0.486$&$0.464/0.411$&$0.419/0.329$&$0.374/0.307$ \\
$0.2$&$512$&$0.620/0.614$&$0.609/0.557$&$0.547/0.455$&$0.476/0.406$&$0.446/0.377$ \\
$0.3$&$512$&$0.888/0.880$&$0.894/0.848$&$0.851/0.804$&$0.823/0.736$&$0.778/0.708$ \\
\noalign{\smallskip}
\noalign{\hrule}
\end{tabular}
\end{center}

\newpage

Table 6:
Empirical power of Monte-Carlo tests $\hat D_m$/$Q_m$ when a first-order autoregressive
model is fit to various indicated autoregressive-moving average models denoted by
models 1--12 in Pe\v{n}a and Rodriguez (2002, Table 3).
The series length was $n=100$ and 1000 simulations were used for each test.

\begin{center}
\begin{tabular}{cccc}
\noalign{\smallskip}
\noalign{\hrule}
\noalign{\smallskip}
\noalign{\hrule}
\noalign{\smallskip}
 model & $m=10$& $m=20$ \\
\noalign{\smallskip}
\noalign{\hrule}
\noalign{\smallskip}
$1$&$0.464$/$0.259$&$0.361$/$0.175$ \\
$2$&$0.988$/$0.705$&$0.974$/$0.539$ \\
$3$&$0.993$/$0.757$&$0.990$/$0.574$ \\
$4$&$0.584$/$0.414$&$0.472$/$0.303$ \\
$5$&$0.584$/$0.414$&$0.472$/$0.303$ \\
$6$&$0.798$/$0.490$&$0.702$/$0.361$ \\
$7$&$1.000$/$0.981$&$1.000$/$0.846$ \\
$8$&$1.000$/$0.821$&$0.998$/$0.621$ \\
$9$&$0.246$/$0.169$&$0.196$/$0.130$ \\
$10$&$0.876$/$0.743$&$0.820$/$0.575$ \\
$11$&$0.252$/$0.119$&$0.192$/$0.083$ \\
$12$&$0.989$/$0.664$&$0.979$/$0.472$ \\
\noalign{\smallskip}
\noalign{\hrule}
\end{tabular}
\end{center}

\newpage

Table 7:
Comparison of $p$-values for portmanteau tests.
The Ljung-Box $Q_m$ test and the Monte-Carlo $\hat D_m$ test are compared
for the ARMA$(2,1)$ model fit to the Ninemile time series and the AR$(2)$ model fit
to Series E.

\begin{center}
\begin{tabular}{cccccc}
\noalign{\smallskip}
\noalign{\hrule}
\noalign{\smallskip}
\noalign{\hrule}
\noalign{\smallskip}
&$m$&$20$&$30$&$40$&50 \\
Ninemile&$Q_m$&$0.9$\%&$8.0$\%&$22.3$\%&$32.2$\% \\
&$\hat D_m$&$0.4$\%&$0.5$\%&$1.1$\%&$1.4$\%  \\
\noalign{\smallskip}
\noalign{\hrule}
\noalign{\smallskip}
&$m$&$5$&$10$&$15$&$20$ \\
Series E&$Q_m$&$4.6$\%&$10.0$\%&$9.0$\%&$21.7$\% \\
&$\hat D_m$&$2.3$\%&$4.3$\%&$4.3$\%&$5.3$\% \\
\noalign{\smallskip}
\noalign{\hrule}

\end{tabular}
\end{center}

\newpage

\begin{figure}[h]
\centerline{\epsfig{figure=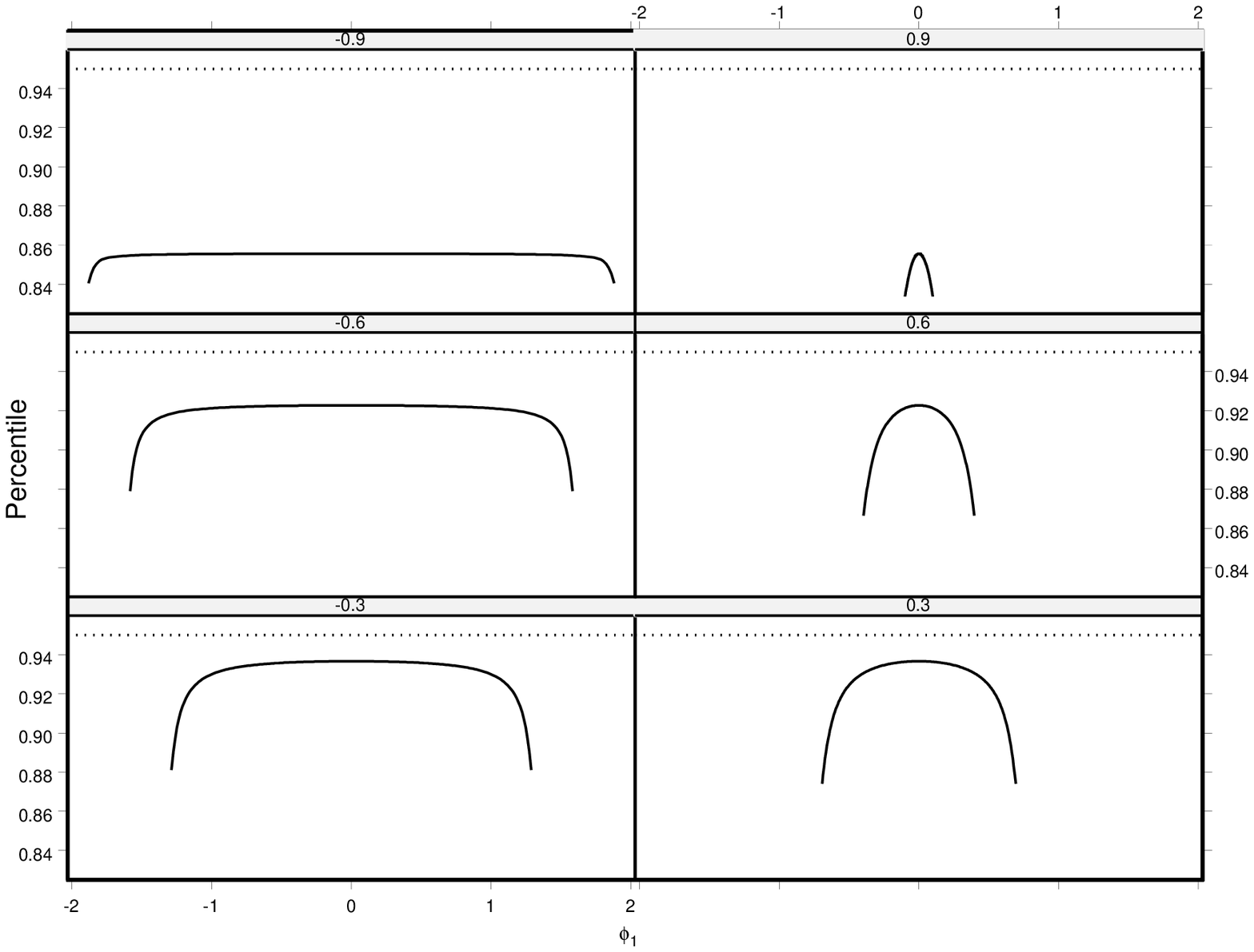,width=5.5in,height=5.5in}}
\caption{
Inaccuracy of the gamma approximation for AR(2) with $m=10$.
$F(Q_\gamma(0.05;\alpha,\beta);\lambda_1,\ldots,\lambda_10)$, where $F$
is the cdf of asymptotic distribution of $\hat D_m$ and $Q_\gamma$ is the quantile
function for the gamma approximation, is plotted.
The plot shows panels corresponding to $\phi_2 = \pm 0.3, \pm 0.6, \pm 0.9$.
The parameter $\phi_1$ extends over the whole admissible
region corresponding to the $\phi_2$.
These plots demonstrate that the gamma approximation is not conservative asymptotically.
When $n$ is large enough, the test based on the gamma approximation will overstate
the true significance level.
}
\end{figure}

\newpage
\begin{figure}
\centerline{\epsfig{figure=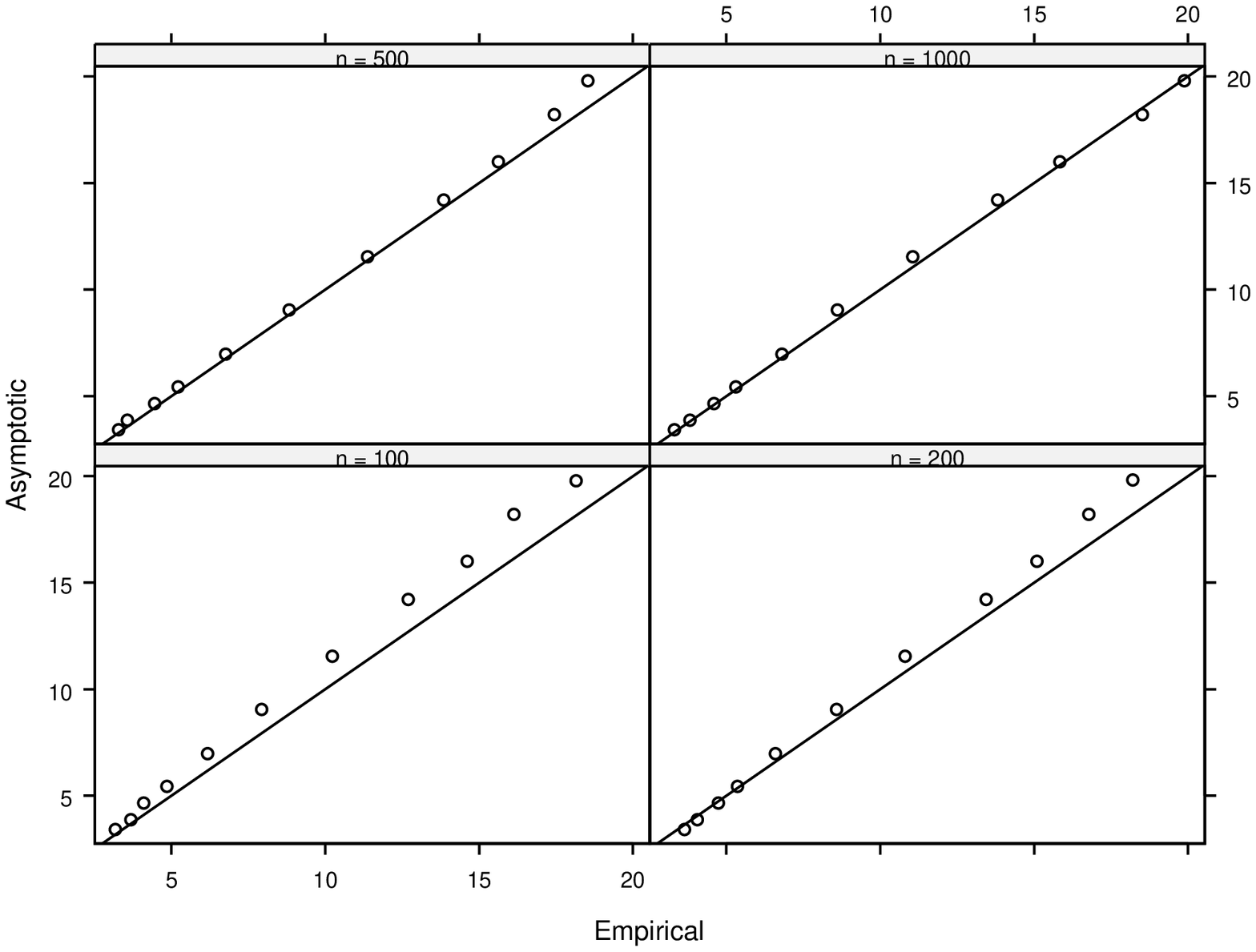,width=6.5in,height=6.5in}}
\caption{Slow convergence of empirical distribution of $\hat D_m$.
The AR$(1)$ with parameter $\phi_1 = 0.5$ was simulated $10^3$ times for series
of lengths $n=100,200,500,1000$ and the empirical distribution of $\hat D_m$ with
$m=20$ was compared with its asymptotic distribution using a QQ plot with quantiles
corresponding to 0.01, 0.02, 0.05, 0.1, 0.2, 0.5, 0.7, 0.9, 0.95, 0.98, 0.99.
}
\end{figure}

\end{document}